\documentclass[10pt]{amsart}
            
\usepackage{amsmath}               
\usepackage{amsfonts}              
\usepackage{amsthm}                
\usepackage[utf8]{inputenc}
\usepackage{leftidx}
\usepackage{tikz-cd}
\usepackage{mathtools}

\newtheorem{tw}{Theorem}[section]
\newtheorem{lem}[tw]{Lemma}
\newtheorem{defi}[tw]{Definition}
\newtheorem{prop}[tw]{Proposition}
\newtheorem{rem}[tw]{Remark}
\newtheorem{cor}[tw]{Corollary}

\newtheorem{ques}[tw]{Question}

\title[Invariance of basic Hodge numbers]{Invariance of basic Hodge numbers under deformations of Sasakian manifolds}
\author{Paweł Raźny}
\address{Institute of Mathematics \\
	Faculty of Mathematics and Computer Science \\
	Jagiellonian University in Cracow
	}
\email{pawel.razny@uj.edu.pl}


\begin{document}

\begin{abstract} We show that the Hodge numbers of Sasakian manifolds are invariant under arbitrary deformations of the Sasakian structure. We also present an upper semi-continuity Theorem for the dimensions of kernels of a smooth family of transversely elliptic operators on manifolds with homologically orientable transversely Riemannian foliations. We use this to prove that the $\partial\bar{\partial}$-lemma and being transversely K\"{a}hler are rigid properties under small deformations of the transversely holomorphic structure which preserve the foliation. Finally, we study an example which shows that this is not the case for arbitrary deformations of the transversely holomorphic foliation.
\end{abstract}
\maketitle

\section{Introduction}
In this short paper we study certain properties of deformations of transversely holomorphic foliations. In \cite{Noz} the authors pose the question whether the basic Hodge numbers of Sasakian manifolds are rigid under arbitrary deformations of Sasakian manifolds. This is motivated by their results on the invariance of such numbers under type I and type II deformations as well as the fact that basic Hodge numbers can be used to distinguish different Sasaki structures on a given manifold. We give a positive answer to the question i.e. we prove the following theorem:
\begin{tw}\label{main} Given a smooth family $\{(M_s,\xi_s,\eta_s,g_s,\phi_s)\}_{s\in [0,1]}$ of compact Sasakian manifolds and fixed integers $p$ and $q$ the function associating to each point $s\in [0,1]$ the basic Hodge number $h^{p,q}_s$ of $(M_s,\xi_s,\eta_s,g_s,\phi_s)$ is constant.
\end{tw}
We split the proof of this result into two theorems which are of independent interest. First we prove Theorem \ref{IfBundle} which states that the basic Hodge numbers are constant for any smooth family (over the interval $[0,1]$) of manifolds with homologically orientable transverse K\"{a}hler foliations for which the spaces of complex-valued basic Harmonic forms constitute a bundle over the interval. Since a family of Sasakian manifolds is in particular a family of homologically orientable transversely K\"{a}hler foliations all that is left to prove is that in this case the spaces of complex-valued basic harmonic forms give in fact a bundle over the interval. This is precisely the content of Theorem \ref{BundleTrue} which allows us to bypass the key difficulty of this and related problems (such as in \cite{Noz}) meaning the fact that the spaces of basic forms over each manifold do not in general form a bundle over the interval. On the way we correct a slight error in \cite{Noz} (see Remark \ref{Correction}). This Theorem strongly relies on the Sasaki structure (and not only on the transverse K\"{a}hler structure) and so the following question remains open:
\begin{ques} Are the basic Hodge numbers rigid under deformations of (homologically orientable) transversely K\"{a}hler foliations on compact manifolds?
\end{ques}
We feel that Theorem \ref{IfBundle} might be helpful in solving this more general problem. Moreover, an answer to this question would have some further use to the theory of $S$-structures which were developed in \cite{B} and are the higher dimensional (meaning the dimension of the characteristic foliation) analogue of Sasakian structures.
\newline\indent In section $4$ we develope some of the Theorems from \cite{KS} for smooth families of transversely elliptic operators on manifolds with TP foliations. We apply this to prove the upper semi-continuity Theorem of the dimensions of kernels of such operators. This in turn is applied to achieve our results in section $5$.
\newline\indent We devote the fifth and sixth section to the study of the behaviour of the basic $\partial\bar{\partial}$-lemma under deformations of transversely holomorphic foliations. We show that if the basic $\partial\bar{\partial}$-lemma holds for a foliated manifold $(M,\mathcal{F})$ then it also holds for appropriately small deformations of the transverse holomorphic structure (provided that we do not deform the foliation itself) as well as a similar rigidity theorem for being transversely K\"{a}hler. The upper-semi continuity theorems for the Bott-Chern and Aeppli cohomology together with the Fr\"{o}licher-type inequality for foliations (which was proven in \cite{My}) allow us to adapt the proofs from \cite{D1} to achieve the main results of this section. In the final section we show that the restriction on deforming the foliation is necessary by studying an example from \cite{Noz,Noz2}.
\section{Preliminaries}
\subsection{Foliations}
We provide a quick review of transverse structures on foliations.
\begin{defi} A codimension q foliation $\mathcal{F}$ on a smooth n-manifold M is given by the following data:
\begin{itemize}
\item An open cover $\mathcal{U}:=\{U_i\}_{i\in I}$ of M.
\item A q-dimensional smooth manifold $T_0$.
\item For each $U_i\in\mathcal{U}$ a submersion $f_i: U_i\rightarrow T_0$ with connected fibers (these fibers are called plaques).
\item For all intersections $U_i\cap U_j\neq\emptyset$ a local diffeomorphism $\gamma_{ij}$ of $T_0$ such that $f_j=\gamma_{ij}\circ f_i$
\end{itemize}
The last condition ensures that plaques glue nicely to form a partition of M consisting of submanifolds of M of codimension q. This partition is called a foliation $\mathcal{F}$ of M and the elements of this partition are called leaves of $\mathcal{F}$.
\end{defi}
We call $T=\coprod\limits_{U_i\in\mathcal{U}}f_i(U_i)$ the transverse manifold of $\mathcal{F}$. The local diffeomorphisms $\gamma_{ij}$ generate a pseudogroup $\Gamma$ of transformations on T (called the holonomy pseudogroup).The space of leaves $M\slash\mathcal{F}$ of the foliation $\mathcal{F}$ can be identified with $T\slash\Gamma$.
\begin{defi}
 A smooth form $\omega$ on M is called basic if for any vector field X tangent to the leaves of $\mathcal{F}$ the following equality holds:
\begin{equation*}
i_X\omega=i_Xd\omega=0
\end{equation*}
Basic 0-forms will be called basic functions henceforth.
\end{defi}
Basic forms are in one to one correspondence with $\Gamma$-invariant smooth forms on T. It is clear that $d\omega$ is basic for any basic form $\omega$. Hence, the set of basic forms of $\mathcal{F}$ (denoted $\Omega^{\bullet}(M\slash\mathcal{F})$) is a subcomplex of the de Rham complex of M. We define the basic cohomology of $\mathcal{F}$ to be the cohomology of this subcomplex and denote it by $H^{\bullet}(M\slash\mathcal{F})$. A transverse structure to $\mathcal{F}$ is a $\Gamma$-invariant structure on T. For example:
\begin{defi}
$\mathcal{F}$ is said to be transversely symplectic if T admits a $\Gamma$-invariant closed 2-form $\omega$ of maximal rank. $\omega$ is then called a transverse symplectic form. As we noted earlier $\omega$ corresponds to a closed basic form of rank q on M (also denoted $\omega$).
\end{defi}
\begin{defi}
$\mathcal{F}$ is said to be transversely holomorphic if T admits a complex structure that makes all the $\gamma_{ij}$ holomorphic. This is equivalent to the existence of an almost complex structure $J$ on the normal bundle $N\mathcal{F}:=TM\slash T\mathcal{F}$ (where $T\mathcal{F}$ is the bundle tangent to the leaves) satisfying:
\begin{itemize}
\item $L_XJ=0$ for any vector field $X$ tangent to the leaves.
\item if $Y_1$ and $Y_2$ are sections of the normal bundle then:
\begin{equation*}
 N_J(Y_1,Y_2):=[JY_1,JY_2]-J[Y_1,JY_2]-J[JY_1,Y_2]+J^2[Y_1,Y_2]=0
 \end{equation*}
where $[$ , $]$ is the bracket induced on the sections of the normal bundle.
\end{itemize}
\end{defi}
\begin{rem}
If $\mathcal{F}$ is transversely holomorphic we have the standard decomposition of the space of complex valued forms $\Omega^{\bullet}({M\slash\mathcal{F},\mathbb{C}})$ into forms of type (p,q) and $d$ decomposes into the sum of operators $\partial$ and $\bar{\partial}$ of order (1,0) and (0,1) respectively. Hence, one can define the Dolbeault double complex $(\Omega^{\bullet,\bullet}({M\slash\mathcal{F},\mathbb{C}}),\partial,\bar{\partial})$, the Fr\"{o}licher spectral sequence and the Dolbeault cohomology as in the manifold case. 
\end{rem}
\begin{defi}
$\mathcal{F}$ is said to be transversely orientable if T is orientable and all the $\gamma_{ij}$ are orientation preserving. This is equivalent to the orientability of $N\mathcal{F}$.
\end{defi}
\begin{defi}
$\mathcal{F}$ is said to be Riemannian if T has a $\Gamma$-invariant Riemannian metric. This is equivalent to the existence of a Riemannian metric g on $N\mathcal{F}$ with $L_Xg=0$ for all vector fields X tangent to the leaves.
\end{defi}
\begin{defi}
$\mathcal{F}$ is said to be transversely parallelizable (TP for short) if there exist q linearly independent $\Gamma$-invariant vector fields.
\end{defi}
Regarding TP foliation we state the following important result from \cite{F}:
\begin{tw}\label{TP} Given a Riemannian TP foliation $\mathcal{F}$ on a compact manifold $M$ the closures of the leaves of $\mathcal{F}$ are submanifolds as well as fibers of a locally trivial fibration $\pi: M\rightarrow W$ with $W$ a compact manifold. In particular they provide another foliation on $M$ for which the leaf space is a compact manifold.
\end{tw}
\begin{defi}
A foliation is said to be Hermitian if it is both transversely holomorphic and Riemannian.
\end{defi}
Throughout the rest of this chapter $\mathcal{F}$ will denote a transversely orientable Riemannian foliation on a compact manifold M. Under these assumptions we shall construct a scalar product on the space of basic forms following \cite{E1}. We start with the principal $SO(q)$-bundle $p:M^{\#}\rightarrow M$ of orthonormal frames transverse to $\mathcal{F}$. The foliation $\mathcal{F}$ lifts to a transversely parallelizable, Riemannian foliation $\mathcal{F}^{\#}$ on $M^{\#}$ of the same dimension as $\mathcal{F}$. Furthermore, this foliation is $SO(q)$-invariant (i.e. for any element $a\in SO(q)$ and any leaf $L$ of $\mathcal{F^{\#}}$, $a(L)$ is also a leaf of $\mathcal{F^{\#}}$) and the transverse metric can be chosen in such a way that it is invariant with respect to the $SO(q)$-action and the fibers of $p:M^{\#}\rightarrow M$ are of measure $1$. By Theorem \ref{TP} there exists a compact manifold $W$ and a fiber bundle $\pi:M^{\#}\rightarrow W$ with fibers equal to the closures of leaves of $\mathcal{F}^{\#}$ (one can now extend the transverse metric to a Riemannian metric on $M^{\#}$ in such a way that the fibers of this bundle have measure $1$ as well). The manifold $W$ is called the basic manifold of $\mathcal{F}$. The $SO(q)$-action on $M^{\#}$ descends to an $SO(q)$-action on $W$. It is apparent that the $SO(q)$-invariant smooth functions on $W$ and basic functions on M are in one to one correspondence. In particular, for basic k-forms $\alpha$ and $\beta$ the basic function $g_x(\alpha_x,\beta_x)$ induces a $SO(q)$-invariant function $\Phi (\alpha,\beta)(w)$ on $W$ (where $g_x$ is the scalar product induced on $\wedge^k T_x^*M$ by the Riemmanian structure). With this we can define the scalar product on basic forms:
\begin{equation*}
<\alpha,\beta>:=\int_W \Phi(\alpha,\beta)(w)d\mu (w)
\end{equation*}
Where $\mu$ is the measure associated to the metric on W. The transverse $*$-operator can be defined fiberwise on the orthogonal complements of the spaces tangent to the leaves in the standard way. This construction can be repeated for complex valued basic forms on Hermitian foliations. We use this scalar product to define $\delta$ as the operator adjoint to $d$ (i.e. such that $<d\alpha,\beta>=<\alpha,\delta\beta>$ for any forms $\alpha$ and $\beta$).
\begin{defi}
A basic differential operator of order m is a linear map $D:\Omega^{\bullet}(M\slash\mathcal{F})\rightarrow\Omega^{\bullet}(M\slash\mathcal{F})$ such that in local coordinates $(x_1,...,x_p,y_1,...,y_q)$ (where $x_i$ are leaf-wise coordinates and $y_j$ are transverse ones) it has the form:
\begin{equation*}
D=\sum\limits_{|s|\leq m}a_s(y)\frac{\partial^{|s|}}{\partial^{s_1}y_1...\partial^{s_q}y_q}
\end{equation*}
where $a_s$ are matrices of appropriate size with basic functions as coefficients. A basic differential operator is called transversely elliptic if its principal symbol is an isomorphism at all points of $x\in M$ and all non-zero, transverse, cotangent vectors at x.
\end{defi}
Due to the correspondence between basic forms of $\mathcal{F}$ and $\Gamma$-invariant forms on the transverse manifold T, a basic differential operator induces a $\Gamma$-invariant differential operator on T. Furthermore, transverse ellipticity of a basic differential operator is equivalent to the ellipticity of its $\Gamma$-invariant counterpart (this is apparent since the principal symbol is defined pointwise).
\begin{tw}(cf.\cite{E1})
Under the above assumptions the kernel of a transversely elliptic differential operator is finitely dimensional.
\end{tw}
\subsection{Basic Bott-Chern and Aeppli cohomology theories}
Let M be a manifold of dimension $n=p+2q$, endowed with a Hermitian foliation $\mathcal{F}$ of complex codimension q. Recall that a foliation satisfies the basic $\partial\bar{\partial}$-lemma if:
$$Ker(\partial)\cap Im(\bar{\partial})=Ker(\bar{\partial})\cap Im(\partial)=Im(\partial\bar{\partial}).$$
This property is thoroughly studied in the classical case in \cite{D1,C1,Del} and in the foliated case in \cite{My}. Suffice to say that in our case it induces many important cohomological properties found in transversely K\"{a}hler foliations such as the decomposition of the basic cohomology induced by the bigradin and the degeneration of the Fr\"{o}licher spectral sequence on the first page.
Using the basic Dolbeault double complex we can define the basic Bott-Chern cohomology of $\mathcal{F}$:
\begin{eqnarray*}
H^{\bullet,\bullet}_{BC}(M\slash\mathcal{F}):=\frac{Ker(\partial)\cap Ker(\bar{\partial})}{Im(\partial\bar{\partial})}
\end{eqnarray*}
where the operators $\partial$ and $\bar{\partial}$ are defined as the components of order (1,0) and (0,1) of the operator $d$ restricted to the basic forms (as mentioned earlier). Our main goal in this subsection, is to present the decomposition theorem for basic Bott-Chern cohomology. To that purpose, we define the operator:
\begin{equation*}
\Delta_{BC}:=(\partial\bar{\partial})(\partial\bar{\partial})^*+(\partial\bar{\partial})^*(\partial\bar{\partial})+
(\bar{\partial}^*\partial)(\bar{\partial}^*\partial)^*+(\bar{\partial}^*\partial)^*(\bar{\partial}^*\partial)
+\bar{\partial}^*\bar{\partial}+\partial^*\partial
\end{equation*}
where by $\partial^*$ and $\bar{\partial}^*$, we mean the operators adjoint to $\partial$ and $\bar{\partial}$,  with respect to the Hermitian product, defined by the transverse Hermitian structure.
\begin{prop}
The operator $\Delta_{BC}$ is transversely elliptic and self-adjoint.
\end{prop}

\begin{tw}\label{HoBC}(Decomposition of basic Bott-Chern cohomology)
If M is a compact manifold, endowed with a Hermitian foliation $\mathcal{F}$, then we have the following decomposition:
\begin{equation*}
\Omega^{\bullet,\bullet}(M\slash\mathcal{F},\mathbb{C})=Ker(\Delta_{BC})\oplus Im(\partial\bar{\partial})\oplus (Im(\partial^*)+Im(\bar{\partial}^*))
\end{equation*}
In particular,
\begin{equation*}
H^{\bullet,\bullet}_{BC}(M\slash\mathcal{F})\cong Ker(\Delta_{BC})
\end{equation*}
and the dimension of $H^{\bullet,\bullet}_{BC}(M\slash\mathcal{F})$ is finite.
\end{tw}

We also define the basic Aeppli cohomology of $\mathcal{F}$ to be:
\begin{equation*}
H^{\bullet,\bullet}_{A}(M\slash\mathcal{F}):=\frac{Ker(\partial\bar{\partial})}{Im(\partial)+Im(\bar{\partial})}
\end{equation*}
We define a basic differential operator, needed for the decomposition theorem for the basic Aeppli cohomology of $\mathcal{F}$:
\begin{equation*}
\Delta_A:=\partial\partial^*+\bar{\partial}\bar{\partial}^*+
(\partial\bar{\partial})^*(\partial\bar{\partial})+(\partial\bar{\partial})(\partial\bar{\partial})^*+(\bar{\partial}\partial^*)^*(\bar{\partial}\partial^*)+
(\bar{\partial}\partial^*)(\bar{\partial}\partial^*)^*
\end{equation*}
\begin{prop}
$\Delta_A$ is a self-adjoint, transversely elliptic operator.
\end{prop}
\begin{tw}(Decomposition of basic Aeppli cohomology)
Let M be a compact manifold, endowed with a Hermitian foliation $\mathcal{F}$. Then we have the following decomposition:
\begin{equation*}
\Omega^{\bullet,\bullet}(M\slash\mathcal{F},\mathbb{C})=Ker(\Delta_A)\oplus (Im(\partial)+Im(\bar{\partial}))\oplus Im((\partial\bar{\partial})^*)
\end{equation*}
In particular, there is an isomorphism,
\begin{equation*}
H^{\bullet,\bullet}_A(M\slash\mathcal{F})\cong Ker(\Delta_A)
\end{equation*}
and the dimension of $H^{\bullet,\bullet}_A(M\slash\mathcal{F})$ is finite.
\end{tw}
Finally, we give a duality theorem for basic Bott-Chern and Aeppli cohomology. However, for the theorem to work, we need an additional condition on our foliation:
\begin{defi} A  codimension $2q$ foliation $\mathcal{F}$ on M is called homologically orientable if $H^{2q}(M\slash\mathcal{F})=\mathbb{R}$.
\end{defi}
\begin{rem} The above condition guaranties that the following equalities hold for basic r-forms:
\begin{equation*}
\partial^*=(-1)^{q(r+1)+1}*\partial*
\quad
\bar{\partial}^*=(-1)^{q(r+1)+1}*\bar{\partial}*
\end{equation*}
where $*$ is the transverse $*$-operator. For general foliations this does not have to be true (c.f. \cite{M1}, appendix B, example 2.3 and \cite{E1}).
\end{rem}

\begin{cor}
If M is a compact manifold endowed with a Hermitian, homologically orientable foliation $\mathcal{F}$, then the transverse star operator induces an isomorphism:
\begin{equation*}
H^{p,q}_{BC}(M\slash\mathcal{F})\rightarrow H^{n-p,n-q}_A(M\slash\mathcal{F})
\end{equation*} 
\end{cor}
Let us continue with the main results from \cite{My}:
\begin{tw}(Basic Fr\"{o}licher-type inequality)
Let $\mathcal{F}$ be a Hermitian foliation of codimension $q$ on a closed manifold M. Then, for every $k\in\mathbb{N}$, the following inequality holds:
\begin{equation*}
\sum\limits_{p+q=k}(dim_{\mathbb{C}}(H^{p,q}_{BC}(M\slash\mathcal{F}))+dim_{\mathbb{C}}(H^{p,q}_A(M\slash\mathcal{F})))\geq 2dim_{\mathbb{C}}(H^k(M\slash\mathcal{F},\mathbb{C}))
\end{equation*}
Furthermore, the equality holds for every $k\in\mathbb{N}$, iff $\mathcal{F}$ satisfies the $\partial\bar{\partial}$-lemma.
\end{tw}
\subsection{Sasakian Manifolds}
We provide a quick recollection of properties of Sasakian Manifolds used in this paper:
\begin{defi} A Sasakian Manifold $(M,g,\xi,\eta,\phi)$ is a $(2n+1)$-dimensional Manifold $M$ together with a Riemannian metric $g$, a Killing vector field $\xi$ a $1$-form $\eta$, and a $(1,1)$ tensor field $\phi$ satisfying for any point $x\in M$ and $X,Y\in T_xM$:
\begin{eqnarray*}
\eta_x\wedge (d\eta^n)_x\neq 0 & \phi^{2}_x(X)=-X+\eta_x(X)\xi_x & \eta_x(\phi_x(X))=0\\
\eta_x(X)=g_x(\xi_x,X) & (d\eta_x)(X,Y)=g_x(\phi_x X,Y)& g_x(\xi_x,\xi_x)=1\\
&g_x(\phi_x(X),\phi_x(Y))=g_x(X,Y)-\eta_x(X)\eta_x(Y)&
\end{eqnarray*}
and additionally the Nijenhuis tensor $[\phi,\phi]$ satisfies:
$$[\phi,\phi](X,Y)+2d\eta(X,Y)\xi=0$$
for any vector fields $X$ and $Y$.
\end{defi}
It is well known that for the homologically orientable foliation $\mathcal{F}$ induced by $\xi$ these tensors define a transverse K\"{a}hler structure by identifying $N\mathcal{F}$ with $\xi^{\perp}$.
\newline\indent Aside from the abundance of properties contained in the above definition and properties of homologically orientable transversely K\"{a}hler foliations we are going to need the following two results:
\begin{prop} For a Sasakian manifold the standard inner product on forms restricted to $\xi^{\perp}$ induced by $g$ can be written in terms of the basic star operator $*_b$ through the formula:
$$<\alpha,\beta>:=\int_M \eta\wedge\alpha\wedge *_b\beta$$
\end{prop}
\begin{tw}\label{BG} Given an odd dimensional manifold $M$ any two Sasaki structures on $M$ have the same basic Betti numbers.
\end{tw}
The later can be found in \cite{BG} (Theorem $7.4.14$).

\section{Invariance of basic Hodge numbers under deformations of Sasakian manifolds}
We start by reducing the problem to proving that the spaces of complex-valued basic harmonic $k$-forms $\mathcal{H}^k_s$ of $(\mathcal{M}_s,\mathcal{F}_s)$ form a bundle over $[0,1]$.
\begin{tw}\label{IfBundle} Let $\{(M_s,\mathcal{F}_s)\}_{s\in [0,1]}$ be a smooth family of homologically orientable transversely K\"{a}hler foliations on compact manifolds such that $\mathcal{H}^k_s$ forms a smooth family of constant dimension for any $k\in\mathbb{N}$. For a fixed pair of integers $(p,q)$ the function associating to each point $s\in [0,1]$ the basic Hodge number $h^{p,q}_s$ of $(M_s,\mathcal{F}_s)$ is constant.
\begin{proof} Using the fact that the kernels of the operators $\Delta$ and $\Delta_{\bar{\partial}}:=\bar{\partial}\bar{\partial}^*+\bar{\partial}^*\bar{\partial}$ are equal under our assumptions (see \cite{E1}) we get the equality:
$$\mathcal{H}^k_s=\bigoplus\limits_{p+q=k}\mathcal{H}^{p,q}_s,$$
where $\mathcal{H}^{p,q}_s$ denotes the kernel of $(\Delta_{\bar{\partial}})_s$ on forms of type $(p,q)$ which is isomorphic to $H^{p,q}(M_s\slash\mathcal{F}_s)$.
Hence, it is sufficient to restrict our attention to the bundle $\mathcal{H}^k_s$. Consider the action of $J_s$ on basic forms given by:
\begin{equation*}
J_s\alpha(X_1,...,X_k)=\sum\limits_{i=1}^k \alpha(X_1,...,J_sX_i,...,X_k),
\end{equation*}
for any $k$ normal sections $X_1,...,X_k\in\Gamma(N\mathcal{F}_s)$ (see e.g. the Lie algebra action in \cite{Gul} for motivation). The spaces $\mathcal{H}^{p,q}_s$ are precisely the $i(p-q)$-eigenspaces of the restriction of $J_s$ to harmonic basic k-forms (note that this operation restricts to a linear operator on $\mathcal{H}^k_s$ due to the decomposition above). With this we can write:
$$\mathcal{H}^{p,q}_s=Ker(J_s|_{\mathcal{H}^k_s}-i(p-q)Id_{\mathcal{H}^k_s}).$$
Taking any $s_0\in [0,1]$ we know (via a standard rank argument) that we can choose a small neighbourhood $U_{p,q}$ of $s_0$ such that the dimension of $Ker(J_s|_{\mathcal{H}^k_s}-i(p-q)Id_{\mathcal{H}^k_s})$ cannot be greater then the dimension of $Ker(J_{s_0}|_{\mathcal{H}^k_{s_0}}-i(p-q)Id_{\mathcal{H}^k_{s_0}})$ for $s\in U_{p,q}$. On the other hand, by our assumptions the direct sum $\bigoplus\limits_{p+q=k}\mathcal{H}^{p,q}_s$ has constant dimension which implies that the dimension of $\mathcal{H}^{p,q}$ cannot drop on $\bigcap\limits_{p+q=k} U_{p,q}$ (since then the dimension of $\mathcal{H}^{p',q'}$ for some other pair $(p',q')$ with $p'+q'=k$ would have to increase to compensate for the loss). This proves that the basic Hodge numbers $h^{p,q}_s$ are locally constant with respect to $s$ and so they are in fact constant. 
\end{proof}
\end{tw}
As we already mentioned in the introduction the main difficulty of the problem is to work around the fact that basic forms may not constitute a bundle over the interval. The first step of dealing with this problem is to consider transverse $k$-forms (i.e. forms $\alpha$ such that $i_{\xi_s}\alpha=0$) we denote the space of such forms by $\Omega^{T,k}_s$ (a similar approach was proposed in e.g. \cite{E2,Noz}). On such forms it is natural to consider the operator $d_T:=\pi(d)$ where $\pi$ is the projection onto transverse forms given by the Riemannian metric. Its adjoint $\delta_T$ is given by the formula:
$$\delta_T:=(-1)^k\star_b^{-1} d_T \star_{b},$$
Which due to homological orientability coincides on basic forms with the basic coderivative $\delta_b$. This allows us to define the transverse Laplace operator in a fashion similar to \cite{E2,Noz}:
$$\Delta^T:=\mathcal{L}_{\xi}\mathcal{L}_{\xi}-\delta_Td_T-d_T\delta_T,$$
and similarly as in \cite{Noz} we can prove the following lemma:
\begin{lem}
The operator $\Delta^T:\Omega^{k,T}\rightarrow \Omega^{k,T}$ is strongly elliptic and self-adjoint.
\begin{proof}
Around any point $x_0$ take a local coordinate chart $(t,x_1,y_1,...,x_n,y_n)$ where $\xi=\frac{\partial}{\partial t}$ and $(x_1,y_1,...,x_n,y_n)$ are transverse holomorphic coordinates such that $(\frac{\partial}{\partial x_1},\frac{\partial}{\partial y_1},...\frac{\partial}{\partial x_n},\frac{\partial}{\partial y_n})$ are orthonormal over $x_0$ and $\eta=dt+\sum\limits_{i=0}^n x_idy_i$. In such coordinates the principal symbol $\sigma(\delta_{T}d_T+d_T\delta_T)$ coincide with that of the Laplacian $\Delta_b$ on the planes $t=0$ (to see this note that in these coordinates $\pi(dt)=-\sum\limits_{i=0}^n x_idy_i$ and so after writing the operator in local coordinates we see that aside from the parts present in $\Delta_b$ the additional components are either of degree less then $2$ or are a multiple of some $x_i$ and hence in either case do not contribute to the symbol over $x_0$). For $\alpha:=\alpha_0dt+\sum\limits_{i=1}^n \alpha_{2i-1}dx_i+\alpha_2idy_i\in T^*_{x_0}M$ let $\sigma_{\alpha}(\Delta_{T})$ be the symbol of $\Delta^T$ at $\alpha$. The symbol $\sigma_{\alpha}(\frac{\partial^2}{\partial^2 t})=\alpha_0^2Id_{(\Omega^{k,T})_{x_0}}$, while the symbol of $\Delta_b$ is given by $\sigma(\Delta_b)=-(\sum\limits_{i=1}^{2n}\alpha_i^2)Id_{(\Omega^{k,T})_{x_0}}$ (see \cite{V} Lemma 5.18). This shows that the symbol $\sigma_{\alpha}(\Delta^T)=||\alpha||^2Id_{(\Omega^{k,T})_{x_0}}$ and so the operator is in fact strongly elliptic.
\newline\indent Since $\delta_{T}d_T+d_T\delta_T$ is self-adjoint it suffices to prove that $\mathcal{L}_{\xi}$ is skew-symmetric. For $\alpha_1,\alpha_2\in\Omega^{k,T}$ we have:
$$\mathcal{L}_{\xi}(\eta\wedge\alpha_1\wedge*_b\overline{\alpha_2})=\eta\wedge\mathcal{L}_{\xi}(\alpha_1)\wedge*_b\overline{\alpha_2}+\eta\wedge\alpha_1\wedge*_b\mathcal{L}_{\xi}\overline{\alpha_2},$$
since $\mathcal{L}_{\xi}\eta=0$ and $\mathcal{L}_{\xi}*_b=*_b\mathcal{L}_{\xi}$. Hence, we only need to prove that the left hand side integrates to zero over $M$. But we can write it as:
$$di_{\xi}(\eta\wedge\alpha_1\wedge*_b\overline\alpha_2)=d(\alpha_1\wedge*_b\overline\alpha_2),$$
now it suffices to note that the right hand-side is exact and hence integrates to zero.
\end{proof}
\end{lem}
\begin{rem}\label{Correction} In \cite{Noz} it is claimed that the form $d(\alpha_1\wedge*_b\overline\alpha_2)$ is itself zero which would also imply our theorem as well as the corresponding theorem in \cite{Noz}. This however is not true since the proof uses transverse forms and not basic ones. More concretely taking $k$-even one can consider the forms $\alpha_1=(d\eta)^{\frac{k}{2}}$ and $\alpha_2=f(d\eta)^{\frac{k}{2}}$ where $f$ is any function on $M$ which is non constant in the $\xi$ direction. It is apparent that $\mathcal{L}_{\xi}(\eta\wedge\alpha_1\wedge*_b\overline\alpha_2)\neq 0$
\newline\indent One can alternatively prove this by computing the adjoint of $\mathcal{L}_{\xi}$ (treated as an operator on $\Omega^{k}(M,\mathbb{C})$) using the formula $\mathcal{L}_{\xi}=di_{\xi}+i_{\xi}d$ along with the fact that $i_{\xi}^*=\eta\wedge$. By standard Hodge theory one arrives at the formula $\eta\wedge=(-1)^k*^{-1}i_{\xi}*$ and using the fact that $\mathcal{L}_{\xi}*=*\mathcal{L}_{\xi}$ one finds that:
$$\mathcal{L}_{\xi}^*=-\mathcal{L}_{\xi}.$$
Now all that is left to prove is that being an adjoint on $\Omega^k(M,\mathbb{C})$ implies being an adjoint on $\Omega^{k,T}$which readily follows from the equalities:
$$\int_M\eta\wedge\mathcal{L}_{xi}\alpha_1\wedge *_b\overline{\alpha}_2=\int_M\mathcal{L}_{xi}\alpha_1\wedge *\overline{\alpha}_2=-\int_M\alpha_1\wedge *\mathcal{L}_{xi}\overline{\alpha}_2=-\int_M\eta\wedge\alpha_1\wedge *_b\mathcal{L}_{xi}\overline{\alpha}_2$$
\end{rem}
With this we can now finish the proof of Theorem \ref{main} by proving the following result:

\begin{tw}\label{BundleTrue} Let $\{(M_s,\xi_s,\eta_s,g_s,\phi_s)\}_{s\in [0,1]}$ be a smooth family of compact Sasakian manifolds over an interval. Then the spaces $\mathcal{H}^k_s$ of complex-valued basic harmonic $k$-forms on $M_s$ constitute a bundle over $[0,1]$.
\begin{proof} We start by using the results of \cite{KS} in a fashion similar to \cite{Noz} in order to contain our problem in some smooth vector bundle (with fibers of finite dimension). Using the Spectral Theorem for smooth families of strongly elliptic self-ajoint operators (see Theorem $1$ of \cite{KS}) for the family $\Delta^{k,T}_s$ we get a complete system of eigensections $\{e_{sh}\}_{h\in{\mathbb{N}}, s\in[0,1]}$ together with the corresponding eigenvalues $\lambda_h(s)$ which form an ascending sequence in $[0,\infty)$ with a single accumulation point at infinity. Fix a point $s_0\in[0,1]$ and let $k_0$ be the largest number such that for $h\in\{1,...,k_0\}$ we have $\lambda_h(s_0)=0$. Consider the family of vector spaces $\mathcal{E}_s=span\{e_{sh}\text{ }|\text{ }h\in\{1,...,k_0\}\}$. Since the only accumulation point of the sequence $\lambda_{h}(s_0)$ is infinity we can find a small disc around $0$ in $\mathbb{C}$ such that the only eigenvalue of $\Delta^{k,T}_{s_0}$ contained in this disc is zero. Using Theorem $2$ of \cite{KS} we establish that for each $h$ the eigenvalues $\lambda_h(s)$ form a continuous function and hence in a small neighbourhood $U$ of $s_0$ all $s\in U$ are contained in this disc as well. This allows us to conclude by using Theorem 3 of \cite{KS} that $P_{\mathcal{E}_s}(\tilde{e}_{sh})$ for $h\in\{1,...,k_0\}$ form smooth sections of $\Omega^{k,T}$ over a small neighbourhood $U'\subset U$ of $s_0$ which span $\mathcal{E}_s$ (where $P_{\mathcal{E}_s}$ is the projection onto $\mathcal{E}_s$ and $\tilde{e}_{sh}$ are the extensions of $e_{s_0h}$ with the use of some partition of unity over $[0,1]$). Shrinking the neighbourhood is necessary to retain linear independence of $\tilde{e}_{sh}$. Hence, we have shown that $\mathcal{E}_s$ form a bundle over $U'$.
\newline\indent Now we consider the operator $\mathcal{L}_{\xi_s}:\mathcal{E}_s\rightarrow \Omega^{k,T}_s$. Note that $Ker\mathcal{L}_{\xi_{s_0}}=\mathcal{H}^k_{s_0}$. Via a standard rank argument there is a small neighbourhood $U''\subset U'$ of $s_0$ such that $dim(Ker\mathcal{L}_{\xi_{s_0}})\geq dim(Ker\mathcal{L}_{\xi_{s}})$. However, $Ker\mathcal{L}_{\xi_{s}}\supset \mathcal{H}^k_s$ and since $dim(\mathcal{H}^k_s)=dim(\mathcal{H}^k_{s_0})$ (by Theorem \ref{BG}) we have the following:
$$dim(Ker\mathcal{L}_{\xi_{s_0}})\geq dim(Ker\mathcal{L}_{\xi_{s}})\geq dim(\mathcal{H}^k_s)=dim(\mathcal{H}^k_{s_0})=dim(Ker\mathcal{L}_{\xi_{s_0}}).$$
Hence, all of the dimensions above are equal and $Ker\mathcal{L}_{\xi_{s}}= \mathcal{H}^k_s$. But this implies that $\mathcal{H}^k_s$ can be described as a kernel of a morphism of bundles and since its dimension is constant we conclude that it is a bundle (over $U''$). It immediately follows that $\mathcal{H}^k_s$ forms a bundle over $[0,1]$ since it is a family of subspaces of a bundle with local trivializations around any point.
\end{proof}
\end{tw}

\section{Upper-semi continuity of dimensions of kernels of transversely elliptic operators}
We start by proving some of the results from \cite{KS} for smooth families of transversely elliptic selfadjoint operators on manifolds with TP foliations. The key tools here are the corresponding Theorems from \cite{KS} as well as methods and constructions from \cite{E1} .
\begin{tw}\label{KS1} Let $M$ be a compact manifold with a codimension $q$ homologically orientable TP Riemannian foliation and let $D:\Omega^k(M\slash\mathcal{F})\rightarrow \Omega^k(M\slash\mathcal{F})$ be a transversely elliptic operator of even order. Then there exists a complete orthonormal set of eigenfunctions $e_h\in \Omega^k(M\slash\mathcal{F})$ with corresponding real eigenvalues $\lambda_h$. Moreover, we can arrange them in such order that the eigenvalues grow and their only possible accumulation point is infinity.
\begin{proof}
First let us note that if our foliation has a dense leaf then the corresponding basic $k$-forms are a finitely dimensional vector space $V$ so the Theorem is trivially true. For a TP foliation it is known that the leaf closures form a bundle over some manifold $W$. Note that there is a natural one to one correspondence between smooth sections of the bundle with fiber over a point $w\in W$ of the form $V_{w}\oplus\Omega^k_w(W)$ (this is the so called useful bundle of \cite{E3}) and basic forms of $\mathcal{F}$. The operator $D$ induces then a self-adjoint elliptic operator $\tilde{D}$ (via this correspondence) acting on the useful bundle over $W$ (in \cite{E1,E2} it was proven that the spaces $V_w$ form a bundle over $W$ and that the operator $\tilde{D}$ acting on this bundle has the desired properties). Now by applying Theorem 1 from \cite{KS} to $\tilde{D}:\Gamma(V)\oplus\Omega^k(W)\rightarrow \Gamma(V)\oplus\Omega^k(W)$ we get our desired result.

\end{proof}
\end{tw}
In the exact same fashion we can adapt Theorems 2 and 3 from \cite{KS} to this context. Hence, we get the following Theorems:
\begin{tw}\label{KS2} Let $M$ be a compact manifold with a codimension $q$ homologically orientable TP Riemannian foliation and let $D_s:\Omega^k(M\slash\mathcal{F})\rightarrow \Omega^k(M\slash\mathcal{F})$ be a family of transversely elliptic operator of even order. Then the eigenvalues $\lambda_h(s)$ in the previous theorem form continuous functions.
\end{tw}
\begin{tw}\label{KS3} Under the assumptions of Theorem \ref{KS2} we put $\mathcal{E}_s:=span\{e_{sh_i}\text{ }|\text{ }i\in\{1,...,l\}\}$ where $e_{sh_i}$ are the eigenfunctions from Theorem \ref{KS1} for the operator $D_s$ such that the corresponding eigenvalues constitute a set of all the eigenvalues contained in some bounded domain $U$ in $\mathbb{C}$ which has no eigenvalues on its boundary. Then the projections onto $\mathcal{E}_s$ depend smoothly on $s$.
\end{tw}
We take the time to pose the following question:
\begin{ques} Can these theorems be further generalized to arbitrary Riemannian foliations?
\end{ques}
Now we are ready to prove the main theorem of this section:
\begin{tw} Let $M$ be a compact manifold with a codimension $q$ homologically orientable Riemannian foliation and let $D_s:\Omega^k(M\slash\mathcal{F})\rightarrow \Omega^k(M\slash\mathcal{F})$ be a family of transversely elliptic operator of even order $m$. Denote $h(s):=dimKer(D)$. Then $h(s)$ is upper semi-continuous.
\begin{proof} We start by lifting the foliation $\mathcal{F}$ to a foliation $\mathcal{F}^{\#}$ on the total space $M^{\#}$ of the bundle of orthonormal frames transverse to $\mathcal{F}$. As we already mentioned in the preliminary section $\mathcal{F}^{\#}$ is TP. Moreover, we have an action of $G=SO(q)$ on $M^{\#}$ such that there is a natural one to one correspondence between $G$-invariant basic forms on $(M^{\#},\mathcal{F}^{\#})$ and basic forms on $(M,\mathcal{F})$ (see \cite{E1} for details). We can now lift the family of operators $D_s$ to a family $D_s^{\#}$ of operators on $M^{\#}$. However, the members of this family are usually not transversely elliptic. To remedy this we consider the family $D_s'$ defined by the formula:
$$D_s':=D_s^{\#}+(-1)^{\frac{m}{2}}(\sum_{i=1}^{N}\mathcal{L}_{Q_i}\mathcal{L}_{Q_i})^{\frac{m}{2}}.$$
where $Q_1,...,Q_N$ are the fundamental vector fields of the $G$-action on $M^{\#}$.
\newline\indent Noting that $Q_i$ are killing (with respect to the transverse metric) via a similar argument as in remark \ref{Correction} we observe that the operators $\mathcal{L}_{Q_i}\mathcal{L}_{Q_i}$ (and hence $D_s'$) are self-adjoint. More precisely, one can prove this by computing the adjoint of $\mathcal{L}_{Q_i}$ (treated as an operator on $\Omega^{k}(M,\mathbb{C})$) using the formula $\mathcal{L}_{Q_i}=di_{Q_i}+i_{Q_i}d$ and the equality $ \delta=(-1)^{k+1}*_bd*_b$ which is true due to homological orientabilty. By standard Hodge theory one arrives at the formula $i_{Q_i}^{*}=(-1)^k*_b^{-1}i_{Q_i}*_b$ and using the fact that $\mathcal{L}_{Q_i}*_b=*_b\mathcal{L}_{Q_i}$ one finds that:
$$\mathcal{L}_{Q_i}^*=-\mathcal{L}_{Q_i}.$$

Due to the results of \cite{E1} the operators $D_s'$ are also strongly transversely elliptic. Note that $D_s'$ coincides on $G$-invariant forms with the operator defined using the identification of $G$-invariant basic forms on $(M^{\#},\mathcal{F}^{\#})$ and basic forms on $(M,\mathcal{F})$.
\newline\indent  We finish the proof by using our adaptations of theorems \cite{KS} in a fashion similar to \cite{Noz} with respect to the family  $D'_s$. Using Theorem \ref{KS1} for the family $D'_s$ we get a complete system of eigensections $\{e_{sh}\}_{h\in{\mathbb{N}}, s\in[0,1]}$ together with the corresponding eigenvalues $\lambda_h(s)$ which form an ascending sequence in $[0,\infty)$ with (at most) a single accumulation point at infinity. Fix a point $s_0\in[0,1]$ and let $k_0$ be the largest number such that for $h\in\{1,...,k_0\}$ we have $\lambda_h{s_0}=0$. Consider the family of vector spaces $\mathcal{E}_s=span\{e_{sh}\text{ }|\text{ }h\in\{1,...,k_0\}\}$. Since the only accumulation point of the sequence $\lambda_{h}(s_0)$ is infinity we can find a small disc around $0$ in $\mathbb{C}$ such that the only eigenvalue of $D'_{s_0}$ contained in this disc is zero. Using Theorem \ref{KS2} we establish that for each $h$ the eigenvalues $\lambda_h(s)$ form a continuous function and hence in a small neighbourhood $U$ of $s_0$ all $s\in U$ are contained in this disc as well. This allows us to conclude by using Theorem \ref{KS3} that $P_{\mathcal{E}_s}(\tilde{e}_{sh})$ for $h\in\{1,...,k_0\}$ form smooth sections of $\Omega^{k}(M\slash\mathcal{F})$ over a small neighbourhood $U'\subset U$ of $s_0$ which span $\mathcal{E}_s$ (where $P_{\mathcal{E}_s}$ is the projection onto $\mathcal{E}_s$ and $\tilde{e}_{sh}$ are the extensions of $e_{s_0h}$ with the use of some partition of unity over $[0,1]$). Shrinking the neighbourhood is necessary to retain linear independence of $\tilde{e}_{sh}$. Hence, we have shown that $\mathcal{E}_s$ form a bundle over $U'$.
\newline\indent Note that since the $G$-action commutes with $D'$ (by the definition of $D'$) and $\mathcal{E}_s$ is a sum of eigenspaces for each $s$ we have a well defined action of $G$ on the family $\mathcal{E}_s$. Let $\mathcal{E}^G_s$ denote the subspace of $\mathcal{E}_s$ consisting of $G$-invariant forms. Due to the fact that representations of compact Lie groups do not change their isomorphism class under smooth deformations we have that $\mathcal{E}^{G}_s$ form a bundle over $U'$. Finally, note that we have:
$$dim(Ker(D_{s_0}))=dim((Ker (D'_{s_0}))^G)=dim(\mathcal{E}^G_s)\geq dim((Ker (D'_{s}))^G)=dim(Ker(D_{s})),$$
for $s\in U'$. This concludes the proof.

\end{proof}
\end{tw}
\begin{rem} Homological orientability is necessary for the self adjointness of $\mathcal{L}_{Q_i}\mathcal{L}_{Q_i}$ as otherwise a correction term appears in the formula for $\delta$. One could remedy this by taking $\mathcal{L}_{Q_i}\mathcal{L}^*_{Q_i}$ instead but then the operators $D'_s$ do not coincide with $D_s$ on basic forms (via the aforementioned correspondence).
\end{rem}
\begin{rem} The above discussion can be easily adapted to complex valued forms and their bi-gradation. Moreover, this can be done even if the transverse holomorphic structures varies with $s$. To see this note that $\pi:(N\mathcal{F})^*\otimes\mathbb{C}\rightarrow (N^{0,1}\mathcal{F})^*_s$ induces an isomorphism between $(N^{0,1}\mathcal{F})^*_{s_0}$ and $(N^{0,1}\mathcal{F})^*_{s_1}$ which preserves basic forms for $s_1$ sufficiently close to $s_0$.
\end{rem}

\begin{cor}\label{uSemi} Let $(M_s,\mathcal{F}_s)$ be a smooth family of compact manifolds with homologically orientable transversely Hermitian foliations such that $\mathcal{F}_{s_1}=\mathcal{F}_{s_2}$ for $s_1,s_2\in [0,1]$ and denote $h^{p,q}_{BC}(s):=dim(H^{p,q}_{BC}(M_s\slash\mathcal{F}_s))$. Then $h^{p,q}_{BC}(s)$ is upper semi-continuous.
\begin{proof} Consider the family $(\Delta_{BC})_s$ of transversely elliptic differential operators. Then by Theorem \ref{HoBC} we have:
$$dim(Ker((\Delta_{BC})_s)|_{\Omega^{p,q}(M_s\slash\mathcal{F}_s,\mathbb{C})})=dim(H^{p,q}_{BC}(M_s\slash\mathcal{F}_s,\mathbb{C})).$$
Hence, after choosing a point $s_0\in [0,1]$ we see that in a sufficiently small neighbourhood $U$ of $s_0$ the dimension of $H^{p,q}_{BC}(M_s\slash\mathcal{F}_s,\mathbb{C})$ can only drop (since they are described as a kernel of a linear operator).
\end{proof}
\end{cor}
\begin{rem} Similar corollaries analogously follow for Dolbeault and Aeppli cohomology theories. One needs to use then the operators $\Delta_{\bar{\partial}}$ and $\Delta_A$.
\end{rem}
\section{Deformations of the transverse holomorphic structure with fixed foliation}
Throughout this section we assume that $\{J_s\}_{s\in[0,1]}$ is a smooth family of transverse Hermitian structures on a compact homologically orientable foliated manifold $(M,\mathcal{F})$ (such deformations were already considered in \cite{G} under the name $f$-deformations). In this section we will show that if $(M,\mathcal{F},J_{s_0})$ satisfies the $\partial\bar{\partial}$-lemma (resp. admits a transverse K\"{a}hler structure) then there exists a neighbourhood $U$ of $s_0$ such that for $s\in U$ the transversly holomorphicly foliated manifold $(M,\mathcal{F},J_{s})$ satisfies the $\partial\bar{\partial}$-lemma (resp. admits a transverse K\"{a}hler structure). We shall show in the subsequent section that this is not the case when the foliation is deformed as well. We will use the notation $(M_s,\mathcal{F}_s)$ instead of $(M,\mathcal{F},J_{s})$ to point out which transverse holomorphic structure is being considered. With the upper semi-continuity theorem of the previous section the rigidity of basic $\partial\bar{\partial}$-lemma is a simple consequence of the foliated version of the Fr\"{o}licher type inequality.
\begin{tw} Let $(M_s,\mathcal{F}_s)$ be a smooth family of compact manifolds with transversely Hermitian homologically orientable foliations such that $\mathcal{F}_{s_1}=\mathcal{F}_{s_2}$ for $s_1,s_2\in [0,1]$. If $(M_{s_0},\mathcal{F}_{s_0})$ satisfies the $\partial\bar{\partial}$-lemma then there exists a neighbourhood $U$ of $s_0$ such that for $s\in U$ the transversly Hermitian foliated manifold $(M_s,\mathcal{F}_s)$ satisfies the $\partial\bar{\partial}$-lemma.
\begin{proof} Using Corollary \ref{uSemi} and the remark that follows we know that the dimensions of both the Bott-Chern and Aeppli cohomologies can only drop on a sufficiently small neighbourhood $U$ of $s_0$. Since $(M_{s_0},\mathcal{F}_{s_0})$ satisfies the $\partial\bar{\partial}$-lemma we have the equality:
$$\sum\limits_{p+q=k}(dim_{\mathbb{C}}(H^{p,q}_{BC}(M_{s_0}\slash\mathcal{F}_{s_0}))+dim_{\mathbb{C}}(H^{p,q}_A(M_{s_0}\slash\mathcal{F}_{s_0})))= 2dim_{\mathbb{C}}(H^k(M_{s_0}\slash\mathcal{F}_{s_0},\mathbb{C})),$$
while the Fr\"{o}licher-type inequality applied to $(M_s,\mathcal{F}_s)$ for $s\in U$ prevents the dimensions of Bott-Chern and Aeppli cohomologies from dropping (since the foliations $\mathcal{F}_s$ coincide for all $s\in [0,1]$).
\end{proof}
\end{tw}
\begin{tw} Let $(M_s,\mathcal{F}_s)$ be a smooth family of compact manifolds with transversely Hermitian homologically orientable foliations such that $\mathcal{F}_{s_1}=\mathcal{F}_{s_2}$ for $s_1,s_2\in [0,1]$. If $(M_{s_0},\mathcal{F}_{s_0})$ is transversely K\"{a}hler then there exists a neighbourhood $U$ of $s_0$ such that for $s\in U$ the transversely Hermitian foliated manifold $(M_s,\mathcal{F}_s)$ is transversely K\"{a}hler.
\begin{proof}Using the Fr\"{o}licher type inequality and Theorem \ref{uSemi} we can again conclude that the dimensions of $Ker(\Delta_{BC})$ are constant (in some small neighbourhood of $s_0$). Using Theorem \ref{KS3} we know that the projection $\pi^{\#}_s:\Omega^{1,1}(M^{\#}\slash\mathcal{F}^{\#})\rightarrow Ker((\Delta'_{BC})_s)$ depends smoothly on $s$. Hence, by restricting $\pi^{\#}$ to $G$-invariant forms and noting that $(\Delta'_{BC})_s$ preserve basic forms we conclude that the same is true for the projection $\pi_s:\Omega^{1,1}(M\slash\mathcal{F})\rightarrow Ker((\Delta_{BC})_s)$. Put:
$$\omega_s:=\frac{1}{2}(\pi_s\omega_{s_0}+\overline{\pi_s\omega_{s_0}}),$$
where $\omega_{s_0}$ is the transverse k\"{a}hler form on $(M_{s_0},\mathcal{F}_{s_0})$ (note that no collision arises since for $\pi_{s_0}$ the expression on the right is in fact equal to $\omega_{s_0}$). Note that the forms $\omega_s$ are real and closed (since they are in $Ker(\Delta_{BC})$). Invariance under $J_s$ follows from being $(1,1)$-forms. Moreover, for $s$ sufficiently close to $s_0$ these forms are non-degenerate and $\omega_s(J_sv,v)\geq 0$. Hence, the forms $\omega_s$ are in fact K\"{a}hler forms for $(M_s,\mathcal{F}_s)$.
\end{proof}
\end{tw}
We also want to present the following simple corollary of the Theorem \ref{uSemi}:
\begin{cor} Let $(M_s,\mathcal{F}_s)$ be a smooth family of compact manifolds with homologically orientable transversely K\"{a}hler foliations such that $\mathcal{F}_{s_1}=\mathcal{F}_{s_2}$ for $s_1,s_2\in [0,1]$. For fixed integers $p$ and $q$ the function associating to each point $s\in [0,1]$ the basic Hodge number $h^{p,q}_s$ of $(M_s,\mathcal{F}_s)$ is constant.
\begin{proof} For transversely K\"{a}hler foliations we have the equality:
$$\sum\limits_{i+j=k}h^{i,j}=h^k,$$
where $h^k$ denotes the basic Betti numbers. Theorem \ref{uSemi} implies that the numbers $h^{i,j}$ cannot increase and hence for the equality to be preserved they have to remain constant.
\end{proof}
\end{cor}
\section{Example}
We recall an example of a family of transversely Hermitian foliations presented in \cite{Noz}. Let $M:=\mathbb{S}^1\times\mathbb{S}^1\times\mathbb{S}^3$. Define $\xi_1$ to be the vector field tangent to the first circle in the product and $\xi_0$ to be the vector field tangent to the fibers of the Hopf fibration (with total space $\mathbb{S}^3$). Let $\xi_s= (1-s)\xi_0+s\xi_1$ and let $\mathcal{F}_s$ denote the foliation of dimension $1$ defined by $\xi_s$. The transverse Hermitian structure is taken from the leaf space $\mathbb{T}^2\times\mathbb{S}^2$ (resp. transverse manifold $\{*\}\times\mathbb{S}^1\times\mathbb{S}^3$) for $s=0$ (resp. $s\in (0,1]$). Equivalently, one can define the transverse holomorphic structure by specyfying the almost complex structure $J$ since the manifold $M$ is parallelizable. In this case one takes $J$ evaluated on the vector field tangent to the second circle to be the orthogonal complement of $\xi_s$ in the tori which are generated by the vector fields $\xi_0$ and $\xi_1$ (the evaluation on the vector fields complementary to $\xi_1$ in the parallelization of $\mathbb{S}^3$ does not change). 
\begin{rem} Since the leaf space of $\mathcal{F}_0$ is precisely $\mathbb{T}^2\times\mathbb{S}^2$ it is in fact transversely K\"{a}hler. For $s\in(0,1]$ the transverse manifold of this foliation can be taken to be $\{*\}\times\mathbb{S}^1\times\mathbb{S}^3$ which admits no closed non-degenerate $2$-form. Hence, these foliations are not transversely symplectic. This proves that being transversely K\"{a}hler is not a rigid property under small deformations if the foliations are allowed to vary.
\end{rem}
Since $(M_0,\mathcal{F}_0)$ is transversely K\"{a}hler it has to also satisfy the $\partial\bar{\partial}$-lemma. We will show that for $s\in [0,1]\backslash\mathbb{Q}$ this lemma does not hold and so we will disprove rigidity of this property when the foliation is allowed to vary.
Note that $\mathbb{S}^3\times\mathbb{S}^1$ is a Lie group which has a basis of one forms $\{\alpha_1,\alpha_2,\alpha_3,\alpha_4\}$ invariant under the action of this group on itself and such that:
$$d\alpha_1=-2\alpha_2\wedge\alpha_3$$
$$d\alpha_2=2\alpha_1\wedge\alpha_3$$
$$d\alpha_3=-2\alpha_1\wedge\alpha_2$$
$$d\alpha_4=0$$
We define the corresponding basis of $(1,0)$ forms by:
$$\beta_1=\alpha_1+i\alpha_2$$
$$\beta_2=\alpha_3+i\alpha_4$$
From this we can easily compute that:
$$\bar{\partial}\beta_1=i\beta_1\wedge\overline{\beta_2}$$
$$\bar{\partial}\beta_2=-i\beta_1\wedge\overline{\beta_1}$$
$$\bar{\partial}\overline{\beta}_1=i\overline{\beta}_1\wedge\overline{\beta_2}$$
$$\bar{\partial}\overline{\beta}_2=0$$
Note that since all these forms are invariant under the action of $\mathbb{S}^3\times\mathbb{S}^1$ they also satisfy $\mathcal{L}_{\xi_s}\beta_i=0$ and hence they are a basis (over $\mathcal{C}^{\infty}(M_s\slash\mathcal{F}_s)$) of basic forms. Note also that for $s\in [0,1]\backslash\mathbb{Q}$ the basic functions are precisely the functions constant in the directions $\xi_0$ and $\xi_1$ hence they can be canonically identified with the functions on $\mathbb{S}^2\times\mathbb{S}^1$. One can now see that $H^{1,1}_A(M\slash\mathcal{F})\neq 0$ since at the very least the form $\beta_2\wedge\overline{\beta}^{2}$ provides a non-vanishing class in it. By the Fr\"{o}licher type inequality it suffices to prove that $H^2(M_s\slash\mathcal{F}_s,\mathbb{C})=0$. By the Fr\"{o}licher spectral sequence it is sufficient to prove that the second Dolbeault cohomology are zero. It is easily seen that $H^{2,0}_{\bar{\partial}}(M_s\slash\mathcal{F}_s)=H^{0,2}_{\bar{\partial}}(M_s\slash\mathcal{F}_s)=0$ since in degree $(2,0)$ the kernel is trivial and for $(0,2)$ the imagis the entire space of $(0,2)$ basic forms. For degree $(1,1)$ one immediately sees that $\bar{\partial}(f\beta_2\wedge\overline{\beta}_1)$ and $\bar{\partial}(f\beta_2\wedge\overline{\beta}_2)$ never vanish while the other two components are contained in the image of $\bar{\partial}$. Hence, we get that for any neighbourhood $U$ of $0$ there exists an $s\in U$ such that $(M_s,\mathcal{F}_s)$ does not satisfy the $\partial\bar{\partial}$-lemma while $(M_0,\mathcal{F}_0)$ satisfies the $\partial\bar{\partial}$-lemma.
\begin{rem} It is important to note that all the foliations in this family are homologically orientable. This is obvious for $s=0$. For $s\in (0,1]$ the generator of the top basic cohomology is provided by $\alpha_1\wedge\alpha_2\wedge\alpha_3\wedge\alpha_4$.
\end{rem}


\begin{thebibliography}{99}
\bibitem{D1} D. Angella,
\emph{Cohomological aspects in complex non-K\"{a}hler geometry}. Springer (2014)
\bibitem{E3} M. Asaoka, A. El Kacimi Alaoui, S. Hurder, K. Richardson,
\emph{Foliations: Dynamics, Geometry and Topology}. Birkh\"{a}user (2014).
\bibitem{B} D.E. Blair,
\emph{Geometry of manifolds with structural group $\mathcal{U}(n)\times\mathcal{O}(s)$}. J. Differential Geom. Vol 4 (2), 155-167 (1970).
\bibitem{BG} C.P. Boyer, K. Galicki,
\emph{Sasakian Geometry} Oxford Mathematical Monographs. Oxford University Press (2007).
\bibitem{C1} G.R. Cavalcanti,
\emph{The decomposition of forms and cohomology of generalized complex manifolds}. J. Geom. Phys. 57(1), 121–132 (2006)
\bibitem{Del} P. Deligne, Ph.A. Griffiths, J. Morgan, D.P. Sullivan,
\emph{Real homotopy theory of K\"{a}hler manifolds}. Invent. Math. 29(3), 245–274 (1975)
\bibitem{E1} A. El Kacimi-Alaoui,
\emph{Op\'{e}rateurs transversalement elliptiques sur un feuilletage riemannien et applications}. Compositio Mathematica, 73, 57-106 (1990).
\bibitem{E2} A. El Kacimi-Alaoui, G. Hector,
\emph{D\"{e}composition de Hodge basique pour un feuilletage riemannien}. Ann. Inst. Fourier 36, 207-227 (1987).
\bibitem{F} E. F\'{e}dida,
\emph{Sur l'existence des feuilletages de Lie} C.R. Acad. Sci. Paris S\'{e}r. A 278 835-837 (1974)
\bibitem{G} J. Girbau, M. Nicolau,
\emph{On deformations of holomorphic foliations}. Annales de l’institut Fourier  39 (2), p. 417-449 (1989).
\bibitem{Noz} O. Goertsches, H. Nozawa, D.T\"{o}ben,
\emph{Rigidity and vanishing of basic Dolbeault cohomology of Sasakian manifolds}. J. Symplect. Geom. 14(1) (2012).
\bibitem{Gul} M. Gualtieri
\emph{Generalized Complex Geometry}. Ph.D. thesis, Oxford University (2003) Math.DG/0401221
\bibitem{Noz2} H. Nozawa,
\emph{Deformation of Sasakian metrics}. T Am Math Soc 366(5) (2008).
\bibitem{KS} K. Kodaira, D.Spencer,
\emph{On Deformations of Complex Analytic Structures III. Stability Theorems for Complex Structures}. Ann. of Math. (2) 71, 43-76 (1960).
\bibitem{M1} P. Molino,
\emph{Riemannian foliations}.  Birkh\"{a}user, 1986. Translated by G. Cairns
\bibitem{My} P. Ra\'{z}ny,
\emph{The Fr\"{o}licher-type inequalities of foliations}, J. Geom. Phys., 114, 593-606 (2017).
\bibitem{V} C. Voisin, Cambridge Studies in Advanced Mathematics, 76. Cambridge University Press (2007).
\emph{Hodge theory and complex algebraic geometry}.
\end{thebibliography}
\end{document}